\documentclass[11pt]{article}
\usepackage{enumerate}
\usepackage{amssymb,a4wide,latexsym,makeidx,epsfig,fleqn}
\usepackage{amsthm}
\usepackage{amsmath}
\usepackage{enumerate}

\begin{document}
\textwidth 150mm \textheight 225mm
\title{On the second smallest and the largest normalized Laplacian eigenvalues of a graph
\thanks{ Supported by
the National Natural Science Foundation of China (No.11171273)}}
\author{{Xiaoguo Tian, Ligong Wang and Yong Lu}\\
{\small Department of Applied Mathematics, School of Science, Northwestern
Polytechnical University,}\\ {\small  Xi'an, Shaanxi 710072,
People's Republic
of China.} \\{\small E-mails: xiaoguotianwycm@163.com; lgwangmath@163.com; luyong.gougou@163.com}}
\date{}
\maketitle
\begin{center}
\begin{minipage}{120mm}
\vskip 0.3cm
\begin{center}
{\small {\bf Abstract}}
\end{center}
{\small Let $G$ be a simple connected graph with order $n$. Let $\mathcal{L}(G)$ be the normalized Laplacian matrix of $G$. Let $\lambda_{k}(G)$ be the $k$-th smallest normalized Laplacian eigenvalue of $G$. Denote $\rho(A)$ the spectral radius of the matrix $A$. In this paper, we study the behaviors of $\lambda_{2}(G)$ and $\rho(\mathcal{L}(G))$ when the graph is perturbed by three operations.

 }
\end{minipage}
\end{center}

\section{Introduction }
\label{sec:ch6-introduction}

Let $A$ be a matrix with order $n\times n$, and $\rho(A)$ be the spectral radius of $A$. Let $G$ be a simple connected graph. Let $V(G)$ and $E(G)$ be the vertex set and the edge set of $G$, respectively. Its order is $|V(G)|$, and its size is $| E(G)|$. For $v\in V(G)$, let $d(v)$ be the degree of $v$, $N_{G}(v)$ be the set of neighbours of a vertex $v$ in $G$. We use the notation $\emph{I}$ for the identity matrix, $e$ for the vector consisting of all ones, $S_{n}$ for the star of order $n$, $C_{n}$ for the cycle of length $n$, $P_{n}$ for the path of length $n-1$ and $Vol(G)$ for the sum of the degrees of all vertices in $G$. Meanwhile, we use the notation $S(G)$ to denote the subdivision graph of $G$, which is the graph obtained from $G$ by inserting some new vertices to some edges of $G$.

   Let $A(G)$ and $D(G)$ be the adjacency matrix and the diagonal matrix of vertex degrees of $G$, respectively. The Laplacian and normalized Laplacian matrices of $G$ are defined as $L(G)=D(G)-A(G)$ and $\mathcal{L}(G)=D^{-\frac{1}{2}}(G)L(G)D^{-\frac{1}{2}}(G)$, respectively. When only one graph $G$ is under consideration, we sometimes use $A$, $D$, $L$ and $\mathcal{L}$ instead of $A(G)$, $D(G)$, $L(G)$ and $\mathcal{L}(G)$ respectively. It is easy to see that $\mathcal{L}(G)$ is a symmertric positive semidefinite matrix and $D^{\frac{1}{2}}(G)e$ is an eigenvector of $\mathcal{L}(G)$ with eigenvalue 0. Thus, the eigenvalues $\lambda_{i}(G)$ of $\mathcal{L}(G)$ satisfy $$\lambda_{n}(G)\geq\cdots\geq\lambda_{2}(G)\geq\lambda_{1}(G)=0.$$ Some of them maybe repeated according to their multiplicities. $\lambda_{k}(G)$ is the $k$-th smallest normalized Laplacian eigenvalue of $G$. Thus $\rho(\mathcal{L}(G))=\lambda_{n}(G)$. When only one graph is under consideration, we may use $\lambda_{k}$ and $\rho(\mathcal{L})$ instead of $\lambda_{k}(G)$ and $\rho(\mathcal{L}(G))$, respectively.

In terms of $\lambda_{2}(G)$, Chung \cite{Chung} showed that $\lambda_{2}(G)$ is $0$ if and only if $G$ is disconnected. This result is closely related to the second smallest eigenvalue of its Laplacian matrix \cite{Butler}. H.H. Li et al. \cite{HHYZ,HH} studied the behavior of $\lambda_{2}$ when the graph is perturbed by grafting an edge and a pendent path, respectively. Recently, J.X. Li et al. \cite{Jian} studied the behavior of $\lambda_{2}$ when the graph is perturbed by separating an edge. They determined all trees and unicyclic graphs with $\lambda_{2}(G)\geq1-\frac{\sqrt{2}}{2}$. Guo et al. \cite{Guo Li} studied the behavior of $\rho(\mathcal{L})$ when the graph is perturbed by removing pendant edges from one vertex to another. The non-bipartite unicyclic
graph with fixed order and girth which has the largest $\rho(\mathcal{L})$ was also determined.

In this paper, we further study the behaviors of $\lambda_{2}$ and $\rho(\mathcal{L})$ when the graph is perturbed by three operations.
\section{Preliminaries}
\label{sec:ch-sufficient}

 In this section we recall some properties of the eigenvalues and eigenfunctions of the normalized Laplacian matrix of a graph $G$. Let $g$ be a vector such that $g\neq0$. Then we can view $g$ as a function which assigns to each vertex $v$ of $G$ a real value $g(v)$, the coordinate of $g$ according to $v$ (All the vectors in this paper are dealt in this way). By letting $g=D^{1/2}f$, we have$$\frac{g^{T}\mathcal{L}g}{g^{T}g}=\frac{f^{T}D^{1/2}\mathcal{L}D^{1/2}f}{(D^{1/2}f)^{T}D^{1/2}f}=\frac{f^{T}Lf}{f^{T}Df}=\frac{\sum_{uv\in E(G)}(f(u)-f(v))^{2}}{\sum_{v\in V(G)}d(v)(f(v))^{2}}.$$ Thus, we can obtain the following formulas for $\lambda_{2}$ and $\rho(\mathcal{L})$.\begin{equation}\label{eq:c1}\lambda_{2}=\inf_{f\perp De}\frac{f^{T}Lf}{f^{T}Df}=\inf_{f\bot De}\frac{\sum_{uv\in E(G)}(f(u)-f(v))^{2}}{\sum_{v\in V(G)}d(v)(f(v))^{2}}.\end{equation}
\begin{equation}\label{eq:c2}\rho(\mathcal{L})=\sup_{f\perp De}\frac{f^{T}Lf}{f^{T}Df}=\sup_{f\bot De}\frac{\sum_{uv\in E(G)}(f(u)-f(v))^{2}}{\sum_{v\in V(G)}d(v)(f(v))^{2}}.\end{equation}

 A nonzero vector that satisfies equality in (1) or (2) is called a
harmonic eigenfunction associated with $\lambda_{2}(G)$ or $\rho(\mathcal{L}(G))$.

\textbf{Lemma 2.1} \cite{Chung} \emph{Let} $G$ \emph{be a simple connected graph and f be a harmonic eigenfunction associated with} $\lambda_{2}(G)$. \emph{Then for any $v\in V(G)$, we have}

\begin{align*}
\frac{1}{d(v)}\sum_{uv\in E(G)}(f(v)-f(u))=\lambda_{2}(G)f(v).
\end{align*}

From Lemma 2.1, we have the following result.

\textbf{Corollary 2.2} \emph{Let $G$ be a simple connected graph and f be a harmonic eigenfunction associated with $\lambda_{2}(G)$. If $\lambda_{2}(G)=1$, then for any $v\in V(G)$, we have }
\begin{align*}
\sum_{uv\in E(G)}f(u)=0.
\end{align*}

If $f$ is a harmonic eigenfunction associated with $\rho(\mathcal{L}(G))$, the similar results about Lemma 2.1 and Corollary 2.2 are obtained.

From Corollary 2.2, we have the following result.

\textbf{Corollary 2.3} \emph{Let $v$ be the center of the star $S_{n}$, f is a harmonic eigenfunction associated  with $\lambda_{2}(S_{n})$. Then $f(v)=0.$}

\textbf{ Proof} Through a simple calculation, we can obtain $\lambda_{2}(S_{n})=1$. Combining the Corollary 2.2, the result is clear. $\square$

\textbf{Lemma 2.4} \cite{Chung} \emph{ For a graph which is not a complete graph, we have} $\lambda_{2}\leq1.$

Next we will define three operations:

\textbf{Operation I}. $G'$ is obtained by inserting a new vertex $w$ to an edge $uv$ of $G$. That is to say $G'=G-uv+uw+wv$.

\textbf{Operation II}. Let $G_{1}$ and $G_{2}$ be two simple connected graphs, $u\in V(G_{1})$, $v\in V(G_{2})$. Let $G$ be a graph obtained from $G_{1}$ and $G_{2}$ by identifying $u$ with $v$ (see Figure 1).

\textbf{Operation III}. Let $u,v$ be two vertices of the simple connected graph $G$. Suppose that $v_{1},v_{2},\ldots,v_{s}$ $(1\leq s\leq d(v))$ are some vertices of $N_{G}(v)\backslash N_{G}(u)$ and $v_{1},v_{2},\ldots,v_{s}$ are different from $u$.
 Let $G'$ be the graph obtained from $G$ by deleting the edges $vv_{i}$ and adding the edges $uv_{i}$. That is to say $G'=G-vv_{1}-vv_{2}-\cdots-vv_{s}+uv_{1}+uv_{2}+\cdots+uv_{s}$.

\begin{figure}[htbp]
  \centering
  \includegraphics[scale=0.55]{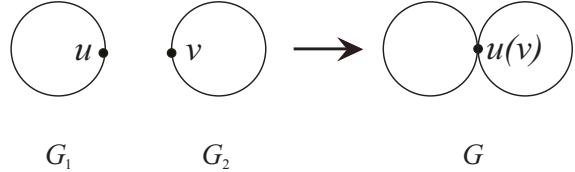}
  \caption{Operation II}
\end{figure}
\section{The effects on the $\lambda_{2}(G)$ of a graph by three operations}
 In this section we study the behavior of $\lambda_{2}$ when the graph is perturbed by three operations.

The following theorem studies the behavior of $\lambda_{2}$ when the graph is perturbed by Operation I.

\textbf{Theorem 3.1}\emph{ Let $G$ be a simple connected graph of order $n$, $uv\in E(G)$ and $G'=G-uv+uw+wv$. Then $\lambda_{2}(G)\geq\lambda_{2}(G')$, and the inequality is strict if $f(u)f(v)\neq0$, where $f$ is a harmonic eigenfunction associated with $\lambda_{2}(G)$.}

\textbf{Proof} Let $V(G)=\{u,v,u_{1},u_{2},\ldots,u_{n-2}\}$ and $V(G')=\{u,v,u_{1},u_{2},\ldots,u_{n-2},$ $w\}$. Let $d(x)$ and $d'(x)$ be the degrees of $x$ in $G$ and $G'$, respectively. Let $D$ and $D'$ be the diagonal degree matrices of $G$ and $G'$, respectively. Let $L$ and $L'$
be the Laplacian matrices of $G$ and $G'$, respectively. Let $e$ and $e'$ be the vectors consisting of all ones, where $e\in R^{n}$ and $e'\in R^{n+1}$. Then $d'(w)=2$ ,$d'(x)=d(x)$, $x\in V(G)$. Since  $f$ is a harmonic eigenfunction associated with $\lambda_{2}(G)$. Then $f\neq0$ and $f\perp De$.
Let us distinguish two cases.

$\textbf{Case 1}$ $f(u)f(v)\leq0$. Let $h$ be a vector such that $h(w)=0$, $h(x)=f(x)$, where $x\in V(G)$. Then
\begin{align*}
\displaystyle
h^{T}D'e'&=\sum_{x\in V(G')}h(x)d'(x)=\sum_{x\in V(G)}h(x)d'(x)+h(w)d'(w)
\end{align*}
\begin{align*}
&=\sum_{x\in V(G)}f(x)d(x)=f^{T}De=0.
\end{align*}
Thus $h\perp D'e'$. Note that $h\neq0$. Then, we have
\begin{align*}
\frac{h^{T}L'h}{h^{T}D'h}\geq\lambda_{2}(G').
\end{align*}
Moreover
\begin{align*}
h^{T}D'h&=\sum_{x\in V(G')}d'(x)h^{2}(x)=\sum_{x\in V(G)}d'(x)h^{2}(x)+d'(w)h^{2}(w)\\
&=\sum_{x\in V(G)}d(x)f^{2}(x)=f^{T}Df,
\end{align*}
and
\begin{align*}
h^{T}L'h=&\sum_{xy\in E(G')}(h(x)-h(y))^{2}\\
&=\sum_{xy\in E(G')\setminus\{uw,wv\}}(h(x)-h(y))^{2}+(h(u)-h(w))^{2}+(h(w)-h(v))^{2}\\
&=\sum_{xy\in E(G)\setminus\{uv\}}(f(x)-f(y))^{2}+f^{2}(u)+f^{2}(v)\\
&=\sum_{xy\in E(G)}(f(x)-f(y))^{2}+2f(u)f(v)=f^{T}Lf+2f(u)f(v)\\
&\leq f^{T}Lf.
\end{align*}
Thus, from Formula \eqref{eq:c1}, we have
\begin{align*}
\lambda_{2}(G)=\frac{f^{T}Lf}{f^{T}Df}\geq\frac{h^{T}L'h}{h^{T}D'h}\geq\lambda_{2}(G').
\end{align*}
If $f(u)f(v)<0$, then $f^{T}Lf>h^{T}L'h$. Thus, $\lambda_{2}(G)>\lambda_{2}(G')$.

$\textbf{Case 2}$ $f(u)f(v)>0$. Let $h$ be a vector such that $h(w)=f(u)$, $h(x)=f(x)$, where $x\in V(G)$. Then
\begin{align*}
h^{T}L'h&=\sum_{xy\in E(G')}(h(x)-h(y))^{2}\\
&=\sum_{xy\in E(G')\setminus\{uw,wv\}}(h(x)-h(y))^{2}+(h(u)-h(w))^{2}+(h(w)-h(v))^{2}\\
&=\sum_{xy\in E(G)\setminus\{uv\}}(f(x)-f(y))^{2}+(f(u)-f(v))^{2}\\
&=\sum_{xy\in E(G)}(f(x)-f(y))^{2}=f^{T}Lf,
\end{align*}
and

\begin{align*}
h^{T}D'e'&=\sum_{x\in V(G')}h(x)d'(x)=\sum_{x\in V(G)}h(x)d'(x)+h(w)d'(w)\\
&=\sum_{x\in V(G)}f(x)d(x)+2f(u)=f^{T}De+2f(u)=2f(u).
\end{align*}
Let $p=h+ce'$, where $c=-\frac{2f(u)}{Vol(G)+2}$. Then
\begin{align*}
p^{T}D'e'=(h+ce')^{T}D'e'=h^{T}D'e'+ce'^{T}D'e'=2f(u)+c(Vol(G)+2)=0.
\end{align*}
Thus $p\perp D'e'$. Note that $p\neq0$. Then, we have
\begin{align*}
\frac{p^{T}L'p}{p^{T}D'p}\geq\lambda_{2}(G').
\end{align*}
It is clear that
\begin{align*}
p^{T}L'p=h^{T}L'h=f^{T}Lf,
\end{align*}
and
\begin{align*}
p^{T}D'p&=\sum_{x\in V(G')}d'(x)p^{2}(x)=\sum_{x\in V(G')}d'(x)(h(x)+c)^{2}\\
&=\sum_{x\in V(G)}d(x)(f(x)+c)^{2}+2(f(u)+c)^{2}\\
&=f^{T}Df+2cf^{T}De+c^{2}Vol(G)+2(f(u)+c)^{2}\\
&=f^{T}Df+\frac{2f^{2}(u)Vol(G)(2+Vol(G))}{(2+Vol(G))^{2}}\\
&>f^{T}Df.
\end{align*}
Thus, from Formula \eqref{eq:c1}, we have
\begin{align*}
\lambda_{2}(G)=\frac{f^{T}Lf}{f^{T}Df}>\frac{p^{T}L'p}{p^{T}D'p}\geq\lambda_{2}(G').
\end{align*}
Combining Cases $1$ and $2$, the result follows.    $\square$

From Theorem 3.1, we have the following result.

\textbf{Corollary 3.2} \emph{Let $G$ be a simple connected graph and $S(G)$ be the subdivision graph of $G$. Then  $\lambda_{2}(G)\geq\lambda_{2}(S(G))$.
}

The following theorem studies the behavior of $\lambda_{2}$ when the graph is perturbed by Operation II.

\textbf{Theorem 3.3} \emph{Let $G_{1}$ and $G_{2}$ be two simple connected graphs of orders $m$ and $n$, respectively. Let $u\in V(G_{1})$ and $v\in V(G_{2})$. Let $G$ be a graph obtained from $G_{1}$ and $G_{2}$ by identifying $u$ with $v$. Then $\lambda_{2}(G)\leq\lambda_{2}(G_{1})$, and the inequality is strict if $f_{1}(u)\neq0$, where $f_{1}$ is a harmonic eigenfunction associated with $\lambda_{2}(G_{1})$.}

\textbf{Proof} Let $V(G_{1})=\{x_{1}, x_{2},\ldots,x_{m-1}, u\}$, $V(G_{2})=\{y_{1}, y_{2},\ldots,y_{n-1}, v\}$, and $V(G)=\{x_{1}, x_{2},\ldots,x_{m-1},u,y_{1},y_{2},\ldots,y_{n-1}\}$. Let $d(x)$, $d_{1}(x)$ and $d_{2}(x)$ be the degree of $x$ in $G$, the degree of $x$ in $G_{1}$, and the degree of $x$ in $G_{2}$, respectively. Let $D$ and $D_{1}$ be the diagonal degree matrices of $G$ and $G_{1}$, respectively. Let $L$ and $L_{1}$
be the Laplacian matrices of $G$ and $G_{1}$, respectively. Let $e$ and $e_{1}$ be the vectors consisting of all ones, where $e\in R^{m+n-1}$ and $e_{1}\in R^{m}$. Then $d(x_{i})=d_{1}(x_{i}), i=1,2,\cdots,m-1$, $d(y_{j})=d_{2}(y_{j}), j=1,2,\ldots,n-1$, and $d(u)=d_{1}(u)+d_{2}(v)$. Since $f_{1}$ is a harmonic eigenfunction associated with $\lambda_{2}(G_{1})$. Then $f_{1}\neq0$ and $f_{1}\perp D_{1}e_{1}$.

Let $f(x)=f_{1}(x),\forall x\in V(G_{1})$, $f(y_{j})=f_{1}(u), j=1,2,\ldots,n-1$. Then we have
\begin{align*}
f^{T}Lf&=\sum_{xy\in E(G)}(f(x)-f(y))^{2}\\
&=\sum_{xy\in E(G_{1})}(f(x)-f(y))^{2}+\sum_{xy\in E(G)\backslash E(G_{1})}(f(x)-f(y))^{2}\\
&=\sum_{xy\in E(G_{1})}(f(x)-f(y))^{2}=\sum_{xy\in E(G_{1})}(f_{1}(x)-f_{1}(y))^{2}\\
&=f_{1}^{T}L_{1}f_{1},
\end{align*}
and
\begin{align*}
f^{T}De&=\sum_{x\in V(G)}d(x)f(x)=\sum_{x\in V(G_{1})}d(x)f(x)+\sum\limits_{j=1}^{n-1}d(y_{j})f(y_{j})\\
&=\sum_{x\in V(G_{1})\backslash\{u\}}d(x)f(x)+d(u)f(u)+\sum\limits_{j=1}^{n-1}d(y_{j})f(y_{j})\\
&=\sum_{x\in V(G_{1})\backslash\{u\}}d_{1}(x)f_{1}(x)+(d_{1}(u)+d_{2}(v))f_{1}(u)+\sum\limits_{j=1}^{n-1}d_{2}(y_{j})f_{1}(u)\\
&=f_{1}^{T}D_{1}e_{1}+d_{2}(v)f_{1}(u)+f_{1}(u)\sum\limits_{j=1}^{n-1}d_{2}(y_{j})\\
&=f_{1}(u)(d_{2}(v)+\sum\limits_{j=1}^{n-1}d_{2}(y_{j}))=f_{1}(u)Vol(G_{2}).
\end{align*}
Let $h=f+ce$, where $c=-\frac{f_{1}(u)Vol(G_{2})}{Vol(G_{1})+Vol(G_{2})}$, Then
\begin{align*}
h^{T}De=(f+ce)^{T}De=f_{1}(u)Vol(G_{2})+c(Vol(G_{1})+Vol(G_{2}))=0.
\end{align*}
Thus $h\perp De$. Note that $h\neq0$. Then, we have
\begin{align*}
\frac{h^{T}Lh}{h^{T}Dh}\geq\lambda_{2}(G).
\end{align*}
Moreover
\begin{align*}
h^{T}Dh&=(f+ce)^{T}D(f+ce)=f^{T}Df+2cf^{T}De+c^{2}e^{T}De\\
\end{align*}
\begin{align*}&=f^{T}Df+2cf_{1}(u)Vol(G_{2})+c^{2}(Vol(G_{1})+Vol(G_{2}))\\
&=f^{T}Df-\frac{(f_{1}(u)Vol(G_{2}))^{2}}{Vol(G_{1})+Vol(G_{2})},
\end{align*}
and
\begin{align*}
f^{T}Df&=\sum_{x\in V(G)}d(x)f^{2}(x)=\sum_{x\in V(G_{1})}d(x)f^{2}(x)+\sum\limits_{j=1}^{n-1}d(y_{j})f^{2}(y_{j})\\
&=\sum_{x\in V(G_{1})\backslash\{u\}}d(x)f^{2}(x)+d(u)f^{2}(u)+\sum\limits_{j=1}^{n-1}d(y_{j})f^{2}(y_{j})\\
&=\sum_{x\in V(G_{1})\backslash\{u\}}d_{1}(x)f^{2}_{1}(x)+(d_{1}(u)+d_{2}(v))f^{2}_{1}(u)+\sum\limits_{j=1}^{n-1}d(y_{j})f^{2}_{1}(u)\\
&=f_{1}^{T}D_{1}f_{1}+d_{2}(v)f^{2}_{1}(u)+f^{2}_{1}(u)\sum\limits_{j=1}^{n-1}d(y_{j})\\
&=f_{1}^{T}D_{1}f_{1}+f^{2}_{1}(u)Vol(G_{2}).
\end{align*}
From the above equation, we have
\begin{align*}
h^{T}Dh&=f_{1}^{T}D_{1}f_{1}+\frac{f^{2}_{1}(u)Vol(G_{1})Vol(G_{2})}{Vol(G_{1})+Vol(G_{2})}\\
&\geq f_{1}^{T}D_{1}f_{1}>0.
\end{align*}
Thus, from Formula \eqref{eq:c1}, we have
\begin{align*}
\lambda_{2}(G_{1})=\frac{f_{1}^{T}L_{1}f_{1}}{f_{1}^{T}D_{1}f_{1}}=\frac{h^{T}Lh}{f_{1}^{T}D_{1}f_{1}}\geq\frac{h^{T}Lh}{h^{T}Dh}\geq\lambda_{2}(G).
\end{align*}
If $f_{1}(u)\neq0$, then $h^{T}Dh>f_{1}^{T}D_{1}f_{1}$. Thus, $\lambda_{2}(G_{1})>\lambda_{2}(G)$. The result follows.$\square$

From Theorem 3.3, the following is easily obtained.

\textbf{Corollary 3.4} \emph{Let $G_{1}$ and $G_{2}$ be two simple connected graphs, $u\in V(G_{1})$, $v\in V(G_{2})$. Let $G$ be a graph obtained from $G_{1}$ and $G_{2}$ by identifying $u$ with $v$. Then $\lambda_{2}(G)\leq min\{\lambda_{2}(G_{1}),\lambda_{2}(G_{2})\}$.
}

In particular, if $T$ is a tree, then it is clear that $T$ can be obtained from the subtree $T_{1}$ and $T_{2}$ by Operation II. Hence, by Corollary 3.4, the following is immediate.

\textbf{Corollary 3.5} \cite{Jianxi} \emph{Let $T$ be a tree. If $T'$ is a subtree of $T$, then $\lambda_{2}(T)\leq\lambda_{2}(T')$}.

The following theorem studies the behavior of $\lambda_{2}$ when the graph is perturbed by Operation III.

\textbf{Theorem 3.6} \emph{Let $u,v$ be two vertices of the simple connected graph $G$ of order $n$. Suppose that $v_{1},v_{2},\ldots,v_{s}(1\leq s\leq d(v))$ are some vertices of $N_{G}(v)\backslash N_{G}(u)$ and $v_{1},v_{2},\ldots,v_{s}$ are different from $u$.
Let $G'=G-vv_{1}-vv_{2}-\cdots-vv_{s}+uv_{1}+uv_{2}+\cdots+uv_{s}$
  and $f$ be a harmonic eigenfunction associated with $\lambda_{2}(G)$. If  $f(u)=f(v)$, Then $\lambda_{2}(G)\geq\lambda_{2}(G')$.}

\textbf{Proof} Let $d(x)$ and $d'(x)$ be the degree of $x$ in $G$ and the degree of $x$ in $G'$, respectively. Let $D$ and $D'$ be the diagonal degree matrices of $G$ and $G'$, respectively. Let $L$ and $L'$
be the Laplacian matrices of $G$ and $G'$, respectively. Let $e$ be the vector consisting of all ones, where $e\in R^{|G|}$. Then $d'(v)=d(v)-s$ ,$d'(u)=d(u)+s$, $d'(x)=d(x)$, where $x\in V(G)\backslash\{u,v\}$. Since $f$ is a harmonic eigenfunction associated with $\lambda_{2}(G)$. Then $f\neq0$ and $f\perp De$.
Let $f'(u)=f(u),\forall u\in V(G)$. Then $f'(u)=f(u)=f(v)=f'(v)$,
\begin{align*}
f'^{T}L'f'&=\sum_{xy\in E(G')}(f'(x)-f'(y))^{2}\\
&=\sum_{xy\in E(G')\backslash\{uv_{1},uv_{2},\cdots,uv_{s}\}}(f'(x)-f'(y))^{2}+\sum\limits_{j=1}^{s}(f'(u)-f'(v_{j}))^{2}\\
&=\sum_{xy\in E(G)\backslash\{vv_{1},vv_{2},\cdots,vv_{s}\}}(f(x)-f(y))^{2}+\sum\limits_{j=1}^{s}(f(v)-f(v_{j}))^{2}\\
&=\sum_{xy\in E(G)}(f(x)-f(y))^{2}\\
&=f^{T}Lf,
\end{align*}
and
\begin{align*}
f'^{T}D'e&=\sum_{x\in V(G')}d'(x)f'(x)\\
&=\sum_{x\in V(G)\backslash\{u,v\}}d'(x)f'(x)+d'(u)f'(u)+d'(v)f'(v)\\
&=\sum_{x\in V(G)\backslash\{u,v\}}d'(x)f'(x)+f'(u)(d'(u)-s)+f'(v)(d'(v)+s)\\
&=\sum_{x\in V(G)\backslash\{u,v\}}d(x)f(x)+d(u)f(u)+d(v)f(v)\\
&=\sum_{x\in V(G)}d(x)f(x)=f^{T}De=0.
\end{align*}
Thus $f'\bot D'e$. Note that $f'\neq0$. Then we have
\begin{align*}
\frac{f'^{T}L'f'}{f'^{T}D'f'}\geq\lambda_{2}(G').
\end{align*}
It is clear that
\begin{align*}
(f')^{T}D'f'&=\sum_{x\in V(G')}d'(x)(f'(x))^{2}\\
&=\sum_{x\in V(G')\backslash\{u,v\}}d'(x)(f'(x))^{2}+d'(u)(f'(u))^{2}+d'(v)(f'(v))^{2}\\
&=\sum_{x\in V(G)\backslash\{u,v\}}d(x)f^{2}(x)+(d(u)+s)f^{2}(u)+(d(v)-s)f^{2}(v)
\end{align*}
\begin{align*}
&=\sum_{x\in V(G)\backslash\{u,v\}}d(x)f^{2}(x)+d(u)f^{2}(u)+d(v)f^{2}(v)\\
&=\sum_{x\in V(G)}d(x)f^{2}(x)=f^{T}Df.
\end{align*}
Hence, from Formula \eqref{eq:c1}, we have
\begin{align*}
\lambda_{2}(G)=\frac{f^{T}Lf}{f^{T}Df}=\frac{(f')^{T}L'f'}{(f')^{T}D'f'}\geq\lambda_{2}(G').
\end{align*}
The result follows.    $\square$

\begin{figure}[htbp]
  \centering
  \includegraphics[scale=0.5]{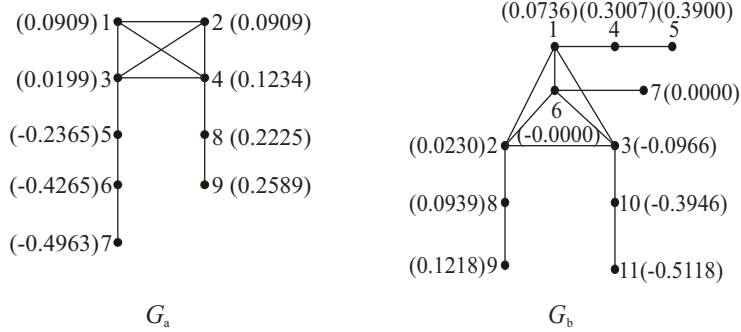}
  \caption{Graph $G_{a}$ and $G_{b}$}
\end{figure}

Figure $2$ shows that the conditoin $f(u)=f(v)$ of Theorem 3.6 is necessary. If $f(u)\neq f(v)$, then the relation between values of $\lambda_{2}(G)$ and $\lambda_{2}(G')$ is not sure. There are the following three cases. For $G_{a}$ and $G_{b}$ in Figure 2 ($G_{b}$ is $G_{3}$ of Theorem 2.4 in \cite{HH}), the natural numbers represent the vertices and the real numbers attached to vertices in each graphs are the valuations by the harmonic eigenfuction associated with $\lambda_{2}(G)$.

\textbf{Case 1} $\lambda_{2}(G)<\lambda_{2}(G')$. Let $G=G_{a}$ and $4$, $3$, $5$ stand for $u$, $v$, $v_{1}$, respectively. It is clear that $f(u)>f(v)$. Denote $G'=G-vv_{1}+uv_{1}$. By direct calculation, we obtain  $\lambda_{2}(G)=0.1408<0.1557=\lambda_{2}(G')$.

\textbf{Case 2} $\lambda_{2}(G)=\lambda_{2}(G')$. Let $G=G_{a}$ and $1$, $3$, $5$ stand for $u$, $v$, $v_{1}$, respectively. It is clear that $f(u)>f(v)$. Denote $G'=G-vv_{1}+uv_{1}$. Because $G$ is isomorphic to $G'$, we obtain  $\lambda_{2}(G)=\lambda_{2}(G')=0.1408$.

\textbf{Case 3} $\lambda_{2}(G)>\lambda_{2}(G')$. Let $G=G_{b}$ and $1$, $2$, $8$ stand for $u$, $v$, $v_{1}$, respectively. It is clear that $f(u)>f(v)$. Denote $G'=G-vv_{1}+uv_{1}$. By direct calculation, we obtain  $\lambda_{2}(G)=0.2290>0.2105=\lambda_{2}(G')$.

From above, we can see that if $f(u)>f(v)$, the relation between values of $\lambda_{2}(G)$ and $\lambda_{2}(G')$ is not sure. Note that if $f$ is a harmonic eigenfunction associated with $\lambda_{2}(G)$, $-f$ is also a harmonic eigenfunction associated with $\lambda_{2}(G)$. When $f(u)<f(v)$, the similar result is obtained.

\section{The effects on the $\rho(\mathcal{L}(G))$ of a graph by three operations}
 In this section we study the behavior of $\rho(\mathcal{L})$ when the graph is perturbed by three operations.

The following theorem studies the behavior of $\rho(\mathcal{L})$ when the graph is perturbed by Operation I.

\textbf{Theorem 4.1}\emph{ Let $G$ be a simple connected graph of order $n$, $uv\in E(G)$ and $G'=G-uv+uw+wv$. Let $f$ be a harmonic eigenfunction associated with $\rho(\mathcal{L}(G))$. If $f(u)f(v)\geq0$ Then $\rho(\mathcal{L}(G))\leq\rho(\mathcal{L}(G'))$, and the inequality is strict if $f(u)f(v)>0$.}

\textbf{Proof} Let $V(G)=\{u,v,u_{1},u_{2},\ldots,u_{n-2}\}$ and $V(G')=\{u,v,u_{1},u_{2},\ldots,u_{n-2},$ $w\}$. Let $d(x)$ and $d'(x)$ be the degrees of $x$ in $G$ and $G'$, respectively. Let $D$ and $D'$ be the diagonal degree matrices of $G$ and $G'$, respectively. Let $L$ and $L'$
be the Laplacian matrices of $G$ and $G'$, respectively. Let $e$ and $e'$ be the vectors consisting of all ones, where $e\in R^{n}$ and $e'\in R^{n+1}$. Then $d'(w)=2$ ,$d'(x)=d(x)$, $x\in V(G)$. Since  $f$ is a harmonic eigenfunction associated with $\rho(\mathcal{L}(G))$. Then $f\neq0$ and $f\perp De$.

Let $h$ be a vector such that $h(w)=0$, $h(x)=f(x)$, where $x\in V(G)$. Then
\begin{align*}
\displaystyle
h^{T}D'e'&=\sum_{x\in V(G')}h(x)d'(x)\\
&=\sum_{x\in V(G)}h(x)d'(x)+h(w)d'(w)\\
&=\sum_{x\in V(G)}f(x)d(x)=f^{T}De=0.
\end{align*}
Thus $h\perp D'e'$. Note that $h\neq0$. Then, we have
\begin{align*}
\frac{h^{T}L'h}{h^{T}D'h}\leq\rho(\mathcal{L}(G')).
\end{align*}
Moreover
\begin{align*}
h^{T}D'h&=\sum_{x\in V(G')}d'(x)h^{2}(x)\\
&=\sum_{x\in V(G)}d'(x)h^{2}(x)+d'(w)h^{2}(w)\\
&=\sum_{x\in V(G)}d(x)f^{2}(x)=f^{T}Df,
\end{align*}
and
\begin{align*}
h^{T}L'h=&\sum_{xy\in E(G')}(h(x)-h(y))^{2}\\
&=\sum_{xy\in E(G')\setminus\{uw,wv\}}(h(x)-h(y))^{2}+(h(u)-h(w))^{2}+(h(w)-h(v))^{2}\\
&=\sum_{xy\in E(G)\setminus\{uv\}}(f(x)-f(y))^{2}+f^{2}(u)+f^{2}(v)
\end{align*}
\begin{align*}
&=\sum_{xy\in E(G)\setminus\{uv\}}(f(x)-f(y))^{2}+(f(u)-f(v))^{2}+2f(u)f(v)\\
&=\sum_{xy\in E(G)}(f(x)-f(y))^{2}+2f(u)f(v)=f^{T}Lf+2f(u)f(v)\\
&\geq f^{T}Lf.
\end{align*}
Thus, from Formula \eqref{eq:c2}, we have
\begin{align*}
\rho(\mathcal{L}(G))=\frac{f^{T}Lf}{f^{T}Df}\leq\frac{h^{T}L'h}{h^{T}D'h}\leq\rho(\mathcal{L}(G')).
\end{align*}
If $f(u)f(v)>0$, then $f^{T}Lf<h^{T}L'h$. Thus, $\rho(\mathcal{L}(G))<\rho(\mathcal{L}(G'))$.     $\square$

It is clear that the proof of Theorem 4.1 is similar to the proof of Case $1$ in Theorem $3.1$.

\begin{figure}[htbp]
  \centering
  \includegraphics[scale=0.5]{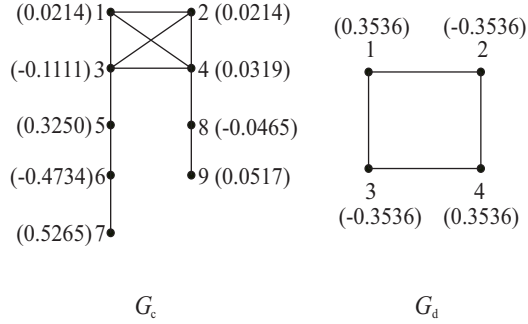}
  \caption{Graph $G_{c}$ and $G_{d}$}
\end{figure}

Figure $3$ shows that the conditoin $f(u)f(v)\geq0$ of Theorem 4.1 is necessary. If $f(u)f(v)<0$, then the relation between values of $\rho(\mathcal{L}(G))$ and $\rho(\mathcal{L}(G'))$ is not sure. There are the following two cases. For $G_{c}$ and $G_{d}$ in Figure 3, the real numbers attached to vertices in each graphs are the valuations by the harmonic eigenfuction associated with $\rho(\mathcal{L}(G))$.

\textbf{Case 1} $\rho(\mathcal{L}(G))<\rho(\mathcal{L}(G'))$. Let $G=G_{c}$ and $5$, $6$ stand for $u$, $v$, respectively. It is clear that $f(u)f(v)<0$. Denote $G'=G-uv+uw+wv$. By direct calculation, we obtain  $\rho(\mathcal{L}(G))=1.8993<1.9382=\rho(\mathcal{L}(G'))$.

\textbf{Case 2} $\rho(\mathcal{L}(G))\geq\rho(\mathcal{L}(G'))$. Let $G=G_{d}$ and $1$, $2$ stand for $u$, $v$, respectively. It is clear that $f(u)f(v)<0$. Denote $G'=G-uv+uw+wv$. By direct calculation, we obtain  $\rho(\mathcal{L}(G))=\rho(\mathcal{L}(C_{4}))=2>1-\cos\frac{4}{5}\pi=\rho(\mathcal{L}(C_{5}))=\rho(\mathcal{L}(G'))$.

The following theorem studies the behavior of $\rho(\mathcal{L})$ when the graph is perturbed by Operation II.

\textbf{Theorem 4.2} \emph{Let $G_{1}$ and $G_{2}$ be two simple connected graphs of orders $m$ and $n$, respectively. Let $u\in V(G_{1})$ and $v\in V(G_{2})$. Let $G$ be a graph obtained from $G_{1}$ and $G_{2}$ by identifying $u$ with $v$. Let $f_{1}$ be a harmonic eigenfunction associated with $\rho(\mathcal{L}(G_{1}))$. If $f_{1}(u)=0$, then $\rho(\mathcal{L}(G))\geq\rho(\mathcal{L}(G_{1}))$. }

\textbf{Proof} Let $V(G_{1})=\{x_{1}, x_{2},\ldots,x_{m-1}, u\}$, $V(G_{2})=\{y_{1}, y_{2},\ldots,y_{n-1}, v\}$, and $V(G)=\{x_{1}, x_{2},\ldots,x_{m-1},u,y_{1},y_{2},\ldots,y_{n-1}\}$. Let $d(x)$, $d_{1}(x)$ and $d_{2}(x)$ be the degree of $x$ in $G$, the degree of $x$ in $G_{1}$, and the degree of $x$ in $G_{2}$, respectively. Let $D$ and $D_{1}$ be the diagonal degree matrices of $G$ and $G_{1}$, respectively. Let $L$ and $L_{1}$
be the Laplacian matrices of $G$ and $G_{1}$, respectively. Let $e$ and $e_{1}$ be the vectors consisting of all ones, where $e\in R^{m+n-1}$ and $e_{1}\in R^{m}$. Then $d(x_{i})=d_{1}(x_{i}), i=1,2,\cdots,m-1$, $d(y_{j})=d_{2}(y_{j}), j=1,2,\ldots,n-1$, and $d(u)=d_{1}(u)+d_{2}(v)$. Since $f_{1}$ is a harmonic eigenfunction associated with $\rho(\mathcal{L}(G_{1}))$. Then $f_{1}\neq0$ and $f_{1}\perp D_{1}e_{1}$.

Let $f(x)=f_{1}(x),\forall x\in V(G_{1})$, $f(y_{j})=0, j=1,2,\ldots,n-1$. Then we have
\begin{align*}
f^{T}Df&=\sum_{x\in V(G)}d(x)f^{2}(x)=\sum_{x\in V(G_{1})}d(x)f^{2}(x)+\sum\limits_{j=1}^{n-1}d(y_{j})f^{2}(y_{j})\\
&=\sum_{x\in V(G_{1})\backslash\{u\}}d(x)f^{2}(x)+d(u)f^{2}(u)\\
&=\sum_{x\in V(G_{1})\backslash\{u\}}d_{1}(x)f^{2}_{1}(x)+(d_{1}(u)+d_{2}(v))f^{2}_{1}(u)\\
&=f_{1}^{T}D_{1}f_{1}+d_{2}(v)f^{2}_{1}(u)=f_{1}^{T}D_{1}f_{1},
\end{align*}

\begin{align*}
f^{T}Lf&=\sum_{xy\in E(G)}(f(x)-f(y))^{2}\\
&=\sum_{xy\in E(G_{1})}(f(x)-f(y))^{2}+\sum_{xy\in E(G)\backslash E(G_{1})}(f(x)-f(y))^{2}\\
&=\sum_{xy\in E(G_{1})}(f_{1}(x)-f_{1}(y))^{2}+\sum_{uy_{j}\in E(G)}(f(u)-f(y_{j}))^{2}\\
&=f_{1}^{T}L_{1}f_{1}+\sum_{uy_{j}\in E(G)}f^{2}(u)\\
&=f_{1}^{T}L_{1}f_{1}+d_{2}(v)f^{2}_{1}(u)=f_{1}^{T}L_{1}f_{1},
\end{align*}

and
\begin{align*}
f^{T}De&=\sum_{x\in V(G)}d(x)f(x)=\sum_{x\in V(G_{1})}d(x)f(x)\\
&=\sum_{x\in V(G_{1})\backslash\{u\}}d(x)f(x)+d(u)f(u)\\
&=\sum_{x\in V(G_{1})\backslash\{u\}}d_{1}(x)f_{1}(x)+(d_{1}(u)+d_{2}(v))f_{1}(u)\\
&=f_{1}^{T}D_{1}e_{1}+d_{2}(v)f_{1}(u)=f_{1}^{T}D_{1}e_{1}=0.
\end{align*}

Thus $f\perp De$. Note that $f\neq0$. Then, we have
\begin{align*}
\frac{f^{T}Lf}{f^{T}Df}\leq\rho(\mathcal{L}(G)).
\end{align*}

Thus, from Formula \eqref{eq:c2}, we have
\begin{align*}
\rho(\mathcal{L}(G_{1}))=\frac{f_{1}^{T}L_{1}f_{1}}{f_{1}^{T}D_{1}f_{1}}=\frac{f^{T}Lf}{f^{T}Df}\leq\rho(\mathcal{L}(G)).
\end{align*}
The result follows.$\square$

It is clear that the proof of Theorem $4.2$ is similar to the proof of Theorem $3.3$.

Figure $3$ shows that the conditoin $f_{1}(u)=0$ of Theorem 4.2 is necessary. If $f_{1}(u)\neq0$, then the relation between values of $\rho(\mathcal{L}(G_{1}))$ and $\rho(\mathcal{L}(G))$ is not sure. There are the following two cases.

\textbf{Case 1} $\rho(\mathcal{L}(G_{1}))<\rho(\mathcal{L}(G))$. Let $G_{1}=G_{c}$, $7$ stand for $u$ and $G_{2}=P_{2}$. It is clear that $f_{1}(u)\neq0$. $G$ is obtained from $G_{1}$ and $G_{2}$ by  Operation II. By direct calculation, we obtain  $\rho(\mathcal{L}(G_{1}))=1.8993<1.9382=\rho(\mathcal{L}(G))$.

\textbf{Case 2} $\rho(\mathcal{L}(G_{1}))\geq\rho(\mathcal{L}(G))$. Let $G_{1}=G_{d}=C_{4}$, $4$ stand for $u$ and $G_{2}=C_{3}$. It is clear that $f_{1}(u)\neq0$. $G$ is obtained from $G_{1}$ and $G_{2}$ by  Operation II. By direct calculation, we obtain  $\rho(\mathcal{L}(G_{1}))=\rho(\mathcal{L}(C_{4}))=2>1.9010=\rho(\mathcal{L}(G))$.

The following theorem studies the behavior of $\rho(\mathcal{L})$ when the graph is perturbed by Operation III.

\textbf{Theorem 4.3} \emph{Let $u,v$ be two vertices of the simple connected graph $G$ of order $n$. Suppose that $v_{1},v_{2},\ldots,v_{s}$ $(1\leq s\leq d(v))$ are some vertices of $N_{G}(v)\backslash N_{G}(u)$ and $v_{1},v_{2},\ldots,v_{s}$ are different from $u$.
Let $G'=G-vv_{1}-vv_{2}-\cdots-vv_{s}+uv_{1}+uv_{2}+\cdots+uv_{s}$,
  and $f$ be a harmonic eigenfunction associated with $\rho(\mathcal{L}(G))$. If $f(u)=f(v)$, then $\rho(\mathcal{L}(G))\leq\rho(\mathcal{L}(G'))$.}

\textbf{Proof} Let $d(x)$ and $d'(x)$ be the degree of $x$ in $G$ and the degree of $x$ in $G'$, respectively. Let $D$ and $D'$ be the diagonal degree matrices of $G$ and $G'$, respectively. Let $L$ and $L'$
be the Laplacian matrices of $G$ and $G'$, respectively. Let $e$ be the vector consisting of all ones, where $e\in R^{|G|}$. Then $d'(v)=d(v)-s$ ,$d'(u)=d(u)+s$, $d'(x)=d(x)$, where $x\in V(G)\backslash\{u,v\}$. Since $f$ is a harmonic eigenfunction associated with $\rho(\mathcal{L}(G))$. Then $f\neq0$ and $f\perp De$.
Let $f'(u)=f(u),\forall u\in V(G)$. Then $f'(u)=f(u)=f(v)=f'(v)$.

Similar to the proof of Theorem $3.6$, we have
\begin{align*}
f'^{T}L'f'=f^{T}Lf, f'^{T}D'f'=f^{T}Df, f'^{T}D'e=f^{T}De=0.
\end{align*}
Thus $f'\bot D'e$. Note that $f'\neq0$. Then we have
\begin{align*}
\frac{f'^{T}L'f'}{f'^{T}D'f'}\leq\rho(\mathcal{L}(G')).
\end{align*}

Thus, from Formula \eqref{eq:c2}, we have
\begin{align*}
\rho(\mathcal{L}(G))=\frac{f^{T}Lf}{f^{T}Df}=\frac{f'^{T}L'f'}{f'^{T}D'f'}\leq\rho(\mathcal{L}(G').
\end{align*}

The result follows.$\square$

Figure $3$ shows that the conditoin $f(u)=f(v)$ of Theorem 4.3 is necessary. If $f(u)\neq f(v)$, then the relation between values of $\rho(\mathcal{L}(G))$ and $\rho(\mathcal{L}(G'))$ is not sure. There are the following three cases.

\textbf{Case 1} $\rho(\mathcal{L}(G))<\rho(\mathcal{L}(G'))$. Let $G=G_{c}$ and $4$, $3$, $5$ stand for $u$, $v$, $v_{1}$, respectively. It is clear that $f(u)>f(v)$. Denote $G'=G-vv_{1}+uv_{1}$. By direct calculation, we obtain  $\rho(\mathcal{L}(G))=1.8993<1.9063=\rho(\mathcal{L}(G'))$.

\textbf{Case 2} $\rho(\mathcal{L}(G))=\rho(\mathcal{L}(G'))$. Let $G=G_{c}$ and $1$, $3$, $5$ stand for $u$, $v$, $v_{1}$, respectively. It is clear that $f(u)>f(v)$. Denote $G'=G-vv_{1}+uv_{1}$. Because $G$ is isomorphic to $G'$, we obtain  $\rho(\mathcal{L}(G))=\rho(\mathcal{L}(G'))=1.8993$.

\textbf{Case 3} $\rho(\mathcal{L}(G))>\rho(\mathcal{L}(G'))$. Let $G=G_{c}$ and $2$, $5$, $6$ stand for $u$, $v$, $v_{1}$, respectively. Considering $g=-f$ as the harmonic eigenfuction associated with $\rho(\mathcal{L}(G))$. It is clear that $g(u)=-0.0214>-0.3250=g(v)$. Denote $G'=G-vv_{1}+uv_{1}$. By direct calculation, we obtain  $\rho(\mathcal{L}(G))=1.8993>1.8243=\rho(\mathcal{L}(G'))$.

From the above, we can see that if $f(u)>f(v)$, then the relation between values of $\rho(\mathcal{L}(G))$ and $\rho(\mathcal{L}(G'))$ is not sure. Note that if $f$ is a harmonic eigenfunction associated with $\rho(\mathcal{L}(G))$, then $-f$ is also a harmonic eigenfunction associated with $\rho(\mathcal{L}(G))$. When $f(u)<f(v)$, the similar result is obtained.

\end{document}